\documentclass[11pt]{amsart}
\parskip=\smallskipamount

\newtheorem{theorem}{Theorem}[section]
\newtheorem{lemma}[theorem]{Lemma}
\newtheorem{corollary}[theorem]{Corollary}
\newtheorem{proposition}[theorem]{Proposition}

\theoremstyle{definition}

\newtheorem{remark}[theorem]{Remark}

\newcommand{\C}{\mathbb{C}}

\newcommand{\N}{\mathbb{N}}

\newcommand{\R}{\mathbb{R}}

\newcommand{\cH}{\mathcal{H}}

\newcommand\wt{\widetilde}
\newcommand\spsh{strongly plurisubharmonic}

\def\di{\partial}
\def\dibar{\overline\partial}
\def\bs{\backslash}

\def\e{\epsilon}

\numberwithin{equation}{section}

\begin{document}
\title[Extending holomorphic mappings]
{Extending holomorphic mappings from subvarieties in Stein manifolds}
\author{Franc Forstneri\v c}
\address{Institute of Mathematics, Physics and Mechanics, 
University of Ljubljana, Jadranska 19, 1000 Ljubljana, Slovenia}
\email{franc.forstneric@fmf.uni-lj.si}
\thanks{Supported by grants P1-0291 and J1-6173, Republic of Slovenia.}

%
%
\date{February 22, 2005} 

\begin{abstract}
Suppose that $Y$ is a complex manifold such that any 
holomorphic map from a compact convex set in a Euclidean space 
$\C^n$ to $Y$ is a uniform limit of entire maps $\C^n\to Y$. 
We prove that a holomorphic map $X_0 \to Y$ from a closed complex 
subvariety $X_0$ in a Stein manifold $X$ admits a holomorphic extension 
$X\to Y$ provided that it admits a continuous extension. 
We then establish the equivalence of four Oka-type properties 
of a complex manifold.
\end{abstract}
\maketitle

\section{Introduction} 
We shall say that a complex manifold $Y$ satisfies the 
{\em convex approximation property} (CAP) if every
holomorphic map $f\colon U \to Y$ from an open set 
$U\subset \C^n$ ($n\in \N$) can be approximated 
uniformly on any compact {\em convex} set $K\subset U$ 
by entire maps $\C^n\to Y$. This Runge type approximation
property was first introduced in \cite{FCAP} 
where it was shown that it implies, and hence is 
equivalent to, the classic {\em Oka property} of $Y$ 
(see corollary \ref{c1} below; in the definition of CAP
we may even  restrict to a certain special class of compact
convex sets). The main result of this paper is that 
CAP also implies the universal extendability of holomorphic 
maps from closed complex subvarieties in Stein manifolds.

\begin{theorem}
\label{t1}
Assume that $Y$ is a complex manifold satisfying {\rm CAP}. 
For any closed complex subvariety $X_0$ in a Stein manifold $X$
and any continuous map $f_0 \colon X\to Y$ such that 
$f_0|_{X_0} \colon X_0\to Y$ is holomorphic
there is a homotopy $f_t\colon X\to Y$ 
$(t\in [0,1])$ which is fixed on $X_0$ such that $f_1$ is 
holomorphic on $X$. The analogous conclusion holds for sections 
of any holomorphic fiber bundle with fiber $Y$
over a Stein manifold. 
\end{theorem}

We shall say that a complex manifold $Y$ satisfying the 
conclusion of theorem \ref{t1} for all data $(X,X_0,f_0)$ 
enjoys the {\em Oka property with interpolation}.
The conclusion of theorem \ref{t1} has been proved
earlier under any of the following assumptions 
(from strongest to weakest):

\begin{itemize}
\item[---] $Y$ is complex homogeneous 
(Grauert \cite{G1}, \cite{G2}, \cite{G3});
\item[---] 
$Y$ admits a {\em dominating spray} (\cite{GOka}, \cite{FP3}); 
such $Y$ is said to be {\em elliptic};
\item[---] $Y$ admits a {\em finite dominating family of sprays}
\cite{Fsubell}; such $Y$ is said to be {\em subelliptic}. 
\end{itemize}
Each of the above conditions implies CAP 
(see the h-Runge theorems proved in \cite{Fsubell},
\cite{FP1}, \cite{GOka}).
The converse implication CAP$\Rightarrow$subellipticity 
is not known in general, and there are cases when 
CAP is known to hold but the existence of a dominating
spray (or a dominating family of sprays) is unclear;
see corollary \ref{c2} below and 
the examples in \cite{FCAP} and \cite{Fflex}.

Our proof shows without any additional work that
the homotopy $\{f_t\}$ in theorem \ref{t1} can be chosen to remain 
holomorphic and uniformly close to $f_0$ on any
compact $\cH(X)$-convex (holomorphically convex) 
subset $K\subset X$ on which the initial map $f_0$ is 
holomorphic (compare with Theorem 1.4 in \cite{FP3}).
However, using an idea of L\'arusson (Theorem 1 in \cite{L}) 
we show that this addition, as well as the jet 
interpolation on complex subvarieties, follows automatically 
(and rather easily) from the Oka property with interpolation. 

\begin{proposition} 
\label{p1}
Assume that a complex manifold $Y$ enjoys the 
Oka property with interpolation
(the conclusion of theorem \ref{t1}). 
Let $d$ be a distance function on $Y$ induced by a Riemannian 
metric. Let $X$ be a Stein manifold, $X_0\subset X$ a closed 
complex subvariety, $K$ a compact $\cH(X)$-convex subset of $X$,
and $f_0\colon X\to Y$ a continuous map which is 
holomorphic in an open set containing $K\cup X_0$. 

For every $s\in \N$ and $\e>0$ there are a neighborhood 
$U\subset X$ of $K\cup X_0$ and a homotopy 
$f_t\colon X\to Y$ $(t\in [0,1])$ such that for every 
$t\in [0,1]$, $f_t$ is holomorphic in $U$, it agrees with 
$f_0$ to order $s$ on $X_0$, and 
$\sup_{x\in K} d(f_t(x),f_0(x)) <\e$.
\end{proposition}

Proposition \ref{p1}, proved in \S 5, easily reduces to the 
case $Y=\C^r$ where it follows from Cartan's theorems A and B. 
The analogous result holds for sections of holomorphic fiber bundles  
over Stein manifolds (remark \ref{p1sections}). 

We shall say that a complex manifold $Y$ satisfying the conclusion 
of proposition \ref{p1} for any $(X,X_0,K,f_0)$ enjoys the 
{\em Oka property with jet interpolation and approximation}. 
Deleting the interpolation (i.e., taking $X_0=\emptyset$)
we get the {\em Oka property with approximation}.
Summarizing the above results we see that these 
ostensibly different Oka-type properties are 
equivalent.

\begin{corollary}
\label{c1}
The following properties of a complex manifold are equivalent:
\begin{enumerate}
\item[(a)] the convex approximation property {\rm (CAP)};
\item[(b)] the Oka property with approximation;
\item[(c)] the Oka property with interpolation;
\item[(d)] the Oka property with jet interpolation and approximation.
\end{enumerate}
\end{corollary}

The analogous equivalences holds for the parametric versions
of these Oka properties (theorem \ref{parametric}).
The conditions in corollary \ref{c1} are implied by ellipticity 
(the existence of a dominating spray; \cite{GOka}, \cite{FP1}) 
and are equivalent to it on Stein manifolds 
(Proposition 1.3 in \cite{FP3}).

The equivalence (a)$\Leftrightarrow$(b) in corollary \ref{c1} 
was proved in \cite{FCAP}, (a)$\Rightarrow$(c) is the content 
of theorem \ref{t1}, (c)$\Rightarrow$(d) is proposition \ref{p1},
and the remaining implications (d)$\Rightarrow$(b) 
and (d)$\Rightarrow$(c) are consequences of definitions.
For (d)$\Rightarrow$(c) we must observe that a holomorphic 
map $X_0\to Y$ from a Stein subvariety $X_0$ of a complex manifold 
$X$ to any complex manifold $Y$ admits a holomorphic extension 
to an open neighborhood of $X_0$ in $X$ (Proposition 1.2 in \cite{FP3}).

In \S 2 we prove an extension of a theorem of Siu \cite{Si} 
on the existence of open Stein neighborhoods of Stein subvarieties. 
This is used to obtain an approxi\-ma\-tion-interpolation 
theorem (\S 3) which is needed in the proof of
theorem \ref{t1} (\S 4). We hope that both these results will be 
of independent interest. Proposition \ref{p1} is proved 
in \S 5. In \S 6 we discuss the analogous equivalences 
for the parametric Oka properties.

%
%
%
%

\section{An extension of a theorem of Siu} 
All complex spaces in this paper are assumed to be reduced and paracompact.
We denote by $\cH(X)$ the algebra of all holomorphic functions 
on a complex space $X$, endowed with the compact-open topology. 
A function (or a map) is said to be {\em holomorphic on a compact set} 
$K$ in $X$ if it is holomorphic in an open set $U\subset X$ containing $K$;
a homotopy $\{f_t\}$ is holomorphic on  $K$ if there is
a neighborhood $U$ of $K$ independent of $t$ such that every 
$f_t$ is holomorphic on $U$. A compact subset $K$ in a Stein space $X$ 
is $\cH(X)$-convex if for any $p\in X\backslash K$ there exists $f\in \cH(X)$ 
with $|f(p)|>\sup_K |f|$. If $K$ is contained in a closed 
complex subvariety $X_0$ of $X$ then $K$
is $\cH(X_0)$-convex if and only if it is $\cH(X)$-convex.
(For the theory of Stein manifolds and Stein spaces we refer 
to \cite{GR} and \cite{Ho}.) 

The following result generalizes a theorem 
of Siu \cite{Si} (1976).

\begin{theorem}
\label{Steinnbds}
Let $X_0$ be a closed Stein subvariety in a complex space $X$.
Assume that $K\subset X$ is a compact set which is 
$\cH(\Omega)$-convex in some open Stein set $\Omega\subset X$ 
containing $K$ and such that $K\cap X_0$ is $\cH(X_0)$-convex.
Then $K\cup X_0$ has a fundamental basis of open 
Stein neighborhoods in $X$.
\end{theorem}

Siu's theorem \cite{Si} corresponds to the case $K=\emptyset$,
the conclusion being that every closed Stein subvariety in a complex 
space admits an open Stein neighborhood. Different proofs of Siu's
theorem were given independently by Col\c toiu \cite{Co}
and Demailly \cite{De} in 1990. If $X$ is Stein and 
$K\subset X$ is $\cH(X)$-convex then $X_0\cup K$ 
admits a basis of Stein neighborhoods which are 
Runge in $X$ (Proposition 2.1 in \cite{CM}).
It seems that under the weaker condition on $K$
in theorem \ref{Steinnbds} the result is new even when $X$ is Stein.
The necessity of $\cH(X_0)$-convexity of  
$K\cap X_0$ is seen by taking
$X=\C^2$, $X_0=\C\times\{0\}$,  and 
$K=\{(z,w) \in\C^2 \colon 1\le |z|\le 2,\ |w|\le 1\}$:
In this case every Stein neighborhood of $K\cup X_0$ 
contains the bidisc $\{(z,w)\colon |z|\le 2,\ |w|\le 1\}$.
(I wish to thank N.\ \O vrelid for this remark.)

\begin{proof} 
We adapt Demailly's proof of Theorem 1 in 
\cite{De}, refering to that paper (or to \cite{Ri})  
for the notion of a strongly plurisubharmonic function
on a complex space with singularities.
(The proof of Col\c toiu \cite{Co} is very similar 
and covers also the more general case when 
$X_0$ is a complete pluripolar set.)
Although we only need the special case 
of theorem \ref{Steinnbds} where $X$ is without singularities, 
the general case does not require any additional effort. 

Let $U\subset X$ be an open set containing $M:=K\cup X_0$.
We shall find an open Stein set $V$ in $X$
with $M\subset V\subset U$. 

By the assumption $K$ is $\cH(\Omega)$-convex in an open 
Stein set $\Omega \subset X$. Hence $K$ has a basis
of open Stein neighborhoods, and replacing $\Omega$ by one 
of them we may assume that $\Omega \subset U$.

Since $K_0:=K\cap X_0$ is assumed to be $\cH(X_0)$-convex, 
it has a compact $\cH(X_0)$-convex neighborhood $K'_0$ 
in $X_0$ which is contained in $\Omega$. Choose a compact 
neighborhood $K'$ of $K$ such that $K'\subset \Omega$ 
and $K'\cap X_0=K'_0$.

Since $K$ is $\cH(\Omega)$-convex, there is a 
smooth strongly plurisubharmonic function $\rho_0$ on 
$\Omega$ such that $\rho_0<0$ on $K$ and $\rho_0>1$ 
on $\Omega\bs K'$ (Theorem 5.1.5 in \cite{Ho}, p.\ 117).
Set $U_c=\{x\in\Omega\colon \rho_0(x)<c\}$.
Fix $c\in (0,1/2)$; then $K\subset U_c\subset U_{2c}\subset K'$.

The restriction $\rho_0|_{X_0\cap \Omega}$ is smooth
strongly plurisubharmonic.  Since the set $K'_0=K'\cap X_0$ 
is assumed to be $\cH(X_0)$-convex, there is a smooth strongly 
plurisubharmonic exhaustion function 
$\rho'_0 \colon X_0\to \R$ which agrees with $\rho_0$ on $K'_0$
and satisfies $\rho'_0 > c$ on $X_0\bs \overline U_c$. 
(To obtain such $\rho'_0$, take a smooth strongly 
plurisubharmonic exhaustion function $\tau\colon X_0\to \R$
such that $\tau<0$ on $K'_0$ and $\tau>1$ on $X_0\bs \Omega$;
also choose a smooth convex increasing function 
$\xi \colon \R\to \R_+$ with $\xi(t)=0$ for $t\le 0$,
and a smooth function $\chi\colon X\to [0,1]$
such that $\chi=1$ on $\{x\in X_0\colon \tau(x)\le 1/2\}$
and $\chi=0$ on $\{x\in X_0\colon \tau(x)\ge 1\}$;
the function $\rho'_0 = \chi\rho_0 + \xi\circ \tau$
satisfies the stated properties provided that $\xi(t)$ is
chosen to grow sufficiently fast for $t>0$.)
Let $\wt \rho_0\colon K'\cup X_0 \to \R$ be defined 
by the conditions $\wt\rho_0|_{K'}=\rho_0|_{K'}$
and $\wt\rho_0|_{X_0}=\rho'_0$.

Choose a smooth convex increasing function $h\colon \R\to \R$ 
satisfying $h(t)\ge t$ for all $t\in \R$, $h(t)=t$ for $t\le c$, 
and $h(t)>t+1$ for $t\ge 2c$. 
The function $\rho_1:=h\circ \wt \rho_0$ is smooth strongly  
plurisubharmonic on $K'\cup X_0$; on the set 
$\overline U_c=\{\rho_0 \le c\}$ we have $\rho_1=\wt \rho_0 =\rho_0$, 
while outside of $U_{2c}$ we have $\rho_1 > \wt\rho_0 +1$. 

By Theorem 4 in \cite{De}, applied to $(\rho_1-1)|_{X_0}$, 
there exists a smooth strongly plurisubharmonic function 
$\rho_2$ in an open neighborhood of $X_0$, satisfying 
$$ 
	\rho_1(x) -1 < \rho_2(x) < \rho_1(x), \qquad x\in X_0.
$$
On a small neighborhood of 
$\overline U_c \cap X_0 =\{x\in X_0\colon \wt \rho_0(x) \le c\}$ 
we have $\rho_2 < \rho_1 = \wt\rho_0$,
while on $X_0\bs U_{2c}$ we have $\rho_2 > \rho_1-1 > \wt\rho_0$. 
It follows that the function
$$
	\rho = \max\{\wt\rho_0,\rho_2\}
$$ 
is well defined and strongly plurisubharmonic in an open 
set $W\subset X$ satisfying $\overline U_c \cup X_0 \subset W \subset U$.
(To see this, observe that the union of the domains of $\wt\rho_0$ and 
$\rho_2$ contains a neighborhood of $\overline U_c \cup X_0$, 
and before running out of the domain of one of these two functions, 
the second function is the larger one and hence takes over.)
After shrinking $W$ around the set $\overline U_c \cup X_0$ 
the function $\rho$ satisfies the following properties:
\begin{itemize}
\item[(i)] $\rho=\wt\rho_0=\rho_0$ on $\overline U_c$
(hence $\rho<0$ on $K$), 
\item[(ii)] $\rho>c$ on $W\bs \overline U_c$, and
\item[(iii)] $\rho=\rho_2$ on $W\bs U_{2c}$.
\end{itemize}
Using the smooth version of the maximum operation as 
in \cite{De} we may also insure that $\rho$ is smooth.
After a further shrinking of $W$ we may assume that 
$\rho$ is a strongly plurisubharmonic exhaustion 
funtion on $\overline W \supset \overline U_c \cup X_0$ 
satisfying $\rho> c$ on $bW$. 
(However, $\rho$ is not an exhaustion function on $W$.) 
The set $L=\{x\in W\colon \rho(x)\le 0\}$ is then compact
and contains $K$ in its interior. 

By Lemma 5 in \cite{De} (see also \cite{Pe})
there is a smooth function $v$ on $X\bs X_0$ 
with a logarithmic pole on $X_0=\{v=-\infty\}$ 
whose Levi form $i\di\dibar v$ is bounded on every 
compact subset of $X$ (such $v$ is said to 
be {\em almost plurisubharmonic}; this notion is 
defined on a complex space by considering 
local ambient extensions as in \cite{Ri}). 
By subtracting a constant from $v$ we may assume 
that $v<0$ on $K$. 

Let $g\colon \R\to \R$ be a convex increasing function 
with $g(t)=t$ for $t\le 0$. For a small $\e>0$ 
(to be specified below) we set 
$$ 
	\wt \rho = \e\, v + g\circ\rho, \quad 
	V=\{x\in W\colon \wt\rho(x)<0\}.
$$
Clearly $\wt \rho|_{X_0}=-\infty$ and hence $X_0\subset V$. 
Furthermore, both summands are negative on $K$ and hence $K\subset V$,
so $V$ is an open neighborhood of $K\cup X_0$ contained in $W$.

To complete the proof we show that $V$ is Stein for a suitable 
choice of $\e$ and $g$. On $L=\{\rho\le 0\}$ we have 
$g\circ\rho=\rho$ which is \spsh; since the Levi form of $v$ 
is bounded on $L$, $\wt \rho$ is 
\spsh\ on $L$ for a sufficiently small $\e>0$. Fix such~$\e$. 
By choosing $g$ to grow sufficiently fast on $(0,+\infty)$  we can 
insure that $\wt\rho$ is \spsh\ on $\overline W$ (since the positive 
Levi form of $g\circ \rho$ can be made sufficiently large 
in order to compensate the bounded negative part of the 
Levi form of $\e v$ on each compact in $\overline W$). 
Furthermore, a sufficiently fast growth of $g$ insures that 
$\wt \rho|_{bW}>0$ and hence $\overline V\subset W$. 

Let $\tau\colon (-\infty,0) \to \R$ be a smooth convex increasing
function such that $\tau(t)=0$ for $t\le -3$ and $\tau(t)=-1/t$
for $t\in (-1,0)$. Then $\psi=\rho + \tau\circ \wt \rho$ is a 
\spsh\ exhaustion function on $V$, and hence $V$ is Stein 
by Narasimhan's theorem \cite{Na}.
\end{proof}

\begin{remark}
\label{additions}
The proof of theorem \ref{Steinnbds} applies also
if $X_0$ is a closed Stein subvariety of $X\bs L$ 
for some compact subset $L\subset {\rm Int}K$,
provided that the set $K\cap X_0$ (which is no longer compact)
is $\cH(X_0)$-convex and $K$ is $\cH(\Omega)$-convex 
in an open Stein neighborhood $\Omega\subset X$. 
For example, if $X_0$ is a complex curve in $X$ 
(possibly with singularities and with some boundary
components in $K$) such that $X_0\bs K$ does not contain any 
irreducible component with compact $X_0$-closure then 
$K\cup X_0$ admits a basis of Stein neigborhoods in $X$. 
An example is obtained by attaching to the solid torus  
$$ K=\{(z,w)\in \C^2\colon 1\le |z|\le 2,\ |w|\le 1\}$$
finitely many punctured discs
$\triangle^*_j=\{(z,b_j)\colon |z|<\frac{3}{2},\ z\ne a_j\}$
where $|a_j|<1$, $|b_j|<1$ and the numbers $b_j$ are distinct ---
the set $K\cup (\cup_{j}\triangle^*_j)$ has a basis of 
Stein neighborhoods (apply theorem \ref{Steinnbds}
with $X=\C^2\bs \cup_j \{(a_j,b_j)\}$).
The analogous conclusion holds if we attach to $K$ finitely many
irreducible complex curves $V_j$ with at least one point 
removed from each $V_j$ in the complement of $K$.
On the other hand, attaching to $K$ a non-punctured disc
one obtains a Hartogs figure without a basis of Stein neighborhoods.
\end{remark}

\section{An approximation--interpolation theorem}
We shall need the following result whose proof relies on 
theorem \ref{Steinnbds}.

\begin{theorem} 
\label{approximation}
Assume that $X$ is a Stein manifold, $X_0\subset X$ a closed 
complex subvariety, $K\subset X$ a compact 
$\cH(X)$-convex set and $U\subset X$ an 
open set containing $K$. Let $Y$ be a complex
manifold with a distance function $d$ induced 
by a Riemannian metric. 
If $f\colon U\cup X_0\to Y$ is a map whose restrictions 
$f|_U$ and $f|_{X_0}$ are holomorphic then for every 
$\e>0$ there exist an open set
$V\subset X$ containing $K\cup X_0$ and a 
holomorphic map $\wt f\colon V\to Y$ such that
$\wt f|_{X_0}=f|_{X_0}$ and 
$\sup_{x\in K} d(\wt f(x),f(x))<\e$.
The analogous result holds for sections
of a holomorphic submersion $Z\to X$.
\end{theorem}

\begin{remark}
Theorem \ref{approximation} holds without any
condition whatsoever on the target manifold $Y$;
however, the domain of $\wt f$ will in general depend 
on $\e$. Our proof gives the analogous result for families 
of maps depending continuously on a parameter in a 
compact Hausdorff space, except that the domains of 
the approximating maps $\wt f$ must be restricted to a
fixed (but arbitrary large) compact in $X$; this suffices for the 
application to the parametric analogue of theorem \ref{t1}.
\end{remark}

\begin{proof}
In the case $Y=\C$ the theorem is well known
and follows from Cartan's theorems A and B: 
The function $f|_{X_0} \colon X_0\to \C$ 
admits a holomorphic extension $\phi\colon X\to \C$;
on $U$ (which we may assume to be Stein and relatively
compact in $X$) can write $f=\phi +\sum_{j=1}^m {g_j h_j}$ 
where the functions $h_j\in \cH(X)$ vanish on $X_0$
and generate the ideal sheaf of $X_0$ over $U$,
and $g_j\in \cH(U)$; approximating $g_j$ uniformly 
on $K$ by $\wt g_j\in \cH(X)$ and taking 
$\wt f=\phi +\sum{\wt g_j h_j}$ completes the proof.

To prove the general case we take $Z=X\times Y$ 
and let $\pi\colon Z\to X$ be the projection $(x,z)\to x$. 
Write $F(x)=(x,f(x))$. The set $Z_0=\{F(x)\in Z \colon x\in X_0\}$
is a closed complex subvariety of $Z$ which is biholomorphic
to $X_0$ (via $F$) and hence is Stein. We may assume that the open 
set $U\subset X$ in theorem \ref{approximation} is Stein and 
$K$ is $\cH(U)$-convex. The set 
$\wt U=\{F(x) \colon x\in U\}$ is a closed Stein
submanifold of $Z|_U=\pi^{-1}(U)$ and hence it has an
open Stein neighborhood $\Omega\subset Z|_U$.

Since $K$ is $\cH(U)$-convex, the set 
$\wt K=\{F(x) \colon x\in K\}$  is $\cH(\wt U)$-convex 
and hence $\cH(\Omega)$-convex. Furthermore, 
$\wt K\cap Z_0=\{F(x)\colon x\in K\cap X_0\}$
is $\cH(Z_0)$-convex since $K\cap X_0$ is 
$\cH(X_0)$-convex. Theorem \ref{Steinnbds} now shows 
that $\wt K \cup Z_0$ has an open Stein neighborhood 
$W\subset Z$. 
Choose a proper holomorphic embedding $\phi\colon W\to\C^r$.
By the special case considered above there is a 
holomorphic map $G\colon X\to \C^r$ such that
$G|_{X_0}=\phi\circ F|_{X_0}$ and $G$ approximates
$\phi\circ F$ on a compact neighborhood of $K$ in $X$.
Let $\iota$ denote a holomorphic retraction 
of a neighborhood of $\phi(W)$ in $\C^r$ onto $\phi(W)$
(Docquier and Grauert \cite{DG}).
The map $\wt F=\phi^{-1}\circ \iota\circ G$ (with values 
in $Z$) is then holomorphic in an open neighborhood 
$V\subset X$ of $K\cup X_0$, it approximates $F$ 
uniformly on $K$ and satisfies $\wt F|_{X_0}=F|_{X_0}$. 
Writing $\wt F(x)=(a(x),\wt f(x)) \in X\times Y$ 
we see that the second component $\wt f$ satisfies the 
conclusion of theorem \ref{Steinnbds}. 

Our proof holds for sections of any holomorphic submersion 
$\pi\colon Z\to X$; when $Z$ does not have a product structure,
the point $\wt F(x)\in Z$ must be projected back to the fiber 
$Z_x=\pi^{-1}(x)$ in order to obtain a 
bona fide section of $\pi\colon Z\to X$. This is accomplished 
by applying to $\wt F(x)$ a holomorphic retraction 
onto $Z_x$ which depends holomorphically on 
the base point $x\in V$ (Lemma 3.4 in \cite{FFourier}).
\end{proof}

%
%
%
%
\section{Proof of theorem \ref{t1}}
The proof which we present requires minimal improvements
of the existing tools, thanks to the new theorem \ref{approximation} 
which provides a local holomorphic extension of a map 
constructed in an inductive step.
Another possibility would be to adapt the proof
of Theorem 1.4 in \cite{FP3} (via Lemma 3.3 and
Proposition 4.2 in \cite{FP3}) to our current situation,
thus proving directly that CAP implies 
the Oka property with jet interpolation and approximation 
(the implication (a)$\Rightarrow$(d) in corollary \ref{c1}).
However, since (c)$\Rightarrow$(d) in corollary \ref{c1}
is quite easy, the latter approach (which seems to require more 
substantial modifications) appears less attractive. 

The main  step is provided by the following.

\begin{proposition} 
\label{pmain}
Let $f_0\colon X\to Y$ and $K,X_0\subset X$ be as in 
theorem \ref{t1}, and let $L\subset X$ be a compact $\cH(X)$-convex 
set with $K\subset {\rm Int}\, L$. For every $\e>0$ there exists
a homotopy $f_s\colon  X\to Y$ $(s\in [0,1])$ 
satisfying the following:
\begin{itemize}
\item[(i)] $f_s|_{X_0}=f_0|_{X_0}$ for every $s\in [0,1]$,  
\item[(ii)] $\sup_{x\in K} d(f_s(x), f_0(x))<\e$ for every $s\in [0,1]$, and 
\item[(iii)] $f_1$ is holomorphic on $L$. 
\end{itemize}
\end{proposition}

Since $X$ is exhausted by a sequence of compact $\cH(X)$-convex 
subsets $K=K_0\subset K_1\subset K_2\subset \ldots$, 
theorem \ref{t1} follows from proposition
\ref{pmain} by an obvious induction.

The remainder of this section is devoted to the proof
of proposition \ref{pmain}. We may assume that 
$L=\{x\in X\colon \rho(x)\le 0\}$ where $\rho\colon X\to\R$
is a smooth strongly plurisubharmonic exhaustion 
function on $X$, $\rho|_K<0$, and $0$ is a regular value of $\rho$.
By theorem \ref{approximation} we may assume that the map
$f_0$ in proposition \ref{pmain} is holomorphic in an open 
set $U\subset X$ containing $K\cup X_0$.

A homotopy satisfying proposition \ref{pmain} is obtained 
by the {\em bumping method} introduced by Grauert in the 
1960's to solve $\dibar$-problems \cite{HL1}. 
To our knowledge this method was first used in the 
Oka-Grauert theory by Henkin and Leiterer in an unpublished
preprint (1986) on which the paper \cite{HL2} is based.
(In these papers the target manifold $Y$ was 
assumed to be complex homogeneous.) Subsequently it has been 
used in most recent works on the subject.
The interpolation requirement presents additional difficulties.
It seems that the first such result (besides those of Grauert) 
is Theorem 1.7 in \cite{FP2}; its proof depends on a geometric 
construction (Lemma 8.4 in \cite{FP2}) which cannot be used 
directly in our situation since it does not insure 
convexity of the bumps. 

The general outline of the bumping method suitable
to our current needs can be found in \S 6 of \cite{FActa} 
and in \S 4 of \cite{FFourier}. The geometric situation 
considered here is precisely the same as in \S 6.5 of \cite{FActa} 
where it was used in the construction of holomorphic submersions 
$X\to\C^q$ with interpolation on a subvariety $X_0\subset X$. 
To avoid repetition as much as possible we shall outline 
the proof (with necessary modifications) and refer to \S 6 
of \cite{FActa} and to \cite{FCAP} for further details where appropriate.  

We recall the geometric setup in \cite{FActa} (p.\ 179);
our current notation is harmonized with that in \cite{FActa}.
The compact set $K' := (K\cup X_0)\cap  \{\rho\le 1\} \subset U$ 
is $\cH(X)$-convex; hence there exists a 
smooth strongly plurisubharmonic exhaustion function 
$\tau \colon X\to \R$ such that $\tau<0$ on $K'$ and $\tau >0$ on $X\bs U$. 
We may assume that $0$ is a regular value of $\tau$ and the hypersurfaces 
$\{\rho=0\}=bL$ and $\{\tau=0\}$ intersect transversely
along the real codimension two submanifold 
$S= \{\rho=0\}\cap \{\tau=0\}$. Hence 
$D_0 :=\{\tau \le 0\} \subset U$ is a smoothly bounded
strongly pseudoconvex domain.
Set $\rho_t=\tau + t(\rho-\tau)=(1-t)\tau + t\rho$ and 
$$
	D_t=\{\rho_t \le 0\} = \{\tau \le t(\tau-\rho) \}, \quad t \in [0,1].
$$
We have $D_0=\{\tau\le 0\}$, $D_1=\{\rho\le 0\}=L$ and $D_0\cap D_1\subset D_t$
for all $t\in [0,1]$. Let 
$$
	\Omega= \{\rho<0,\ \tau>0\} \subset D_1\bs D_0 \quad {\rm and} \quad 
        \Omega'=\{\rho>0,\ \tau<0\} \subset D_0\bs D_1.
$$
As $t$ increases from $0$ to $1$, $D_t \cap L$  increases 
to $D_1=L$ while $D_t\bs L \subset D_0$ decrease to $\emptyset$. 
All hypersurfaces $\{\rho_t=0\}=bD_t$ intersect along 
$S= \{\rho=0\}\cap \{\tau=0\}$. Since 
$d\rho_t=(1-t)d\tau+t d\rho$ and the differentials $d\tau$, $d\rho$
are linearly independent along $S$, $bD_t$ is smooth near $S$. 
Finally, $\rho_t$ is strongly plurisubharmonic and hence $D_t$ is 
strongly pseudoconvex at every smooth boundary point for every $t\in [0,1]$.

We find a homotopy satisfying proposition 
\ref{pmain} by a sequence of finitely many localized steps.
Fix two nearby values of the parameter, say $0\le t_0<t_1 \le 1$.  
Consider first the case that all boundaries $bD_t$ 
for $t\in [t_0,t_1]$ are smooth. Then $D_{t_1}$ 
is obtained from $D_{t_0}$ by attaching finitely many small convex bumps 
of the type described in \cite{FActa} (disjoint from $X_0$)
and intersecting the union with the set $D_{t_1}$. 
Here `small convex' refers to a suitable local holomorphic 
coordinate system on $X$. 
At every step we have a map $X\to Y$ which is holomorphic 
on a certain compact set $A$ obtained from $D_{t_0}$ by the 
previous attachments of bumps; the goal of the step is to obtain 
a new map which is holomorphic on $A\cup B$
(where $B$ is the next convex bump), approximates the given map 
uniformly on $A$ and agrees with it on $X_0$.
The pair $(A,B)$ is chosen such that $C=A\cap B$ is convex 
in some local holomorphic coordinates in a neighborhood 
of $B$ in $X$, $A\cup B$ admits a Stein neighborhood basis in $X$, 
and $\overline {A\backslash B} \cap \overline {B\backslash A}=\emptyset$. 
(To prove an effective version of theorem \ref{t1} in which the CAP axiom 
is used only for maps from Euclidean spaces of dimension 
$\le \dim X+\dim Y$ to $Y$, we must assume in addition that 
$C$ is Runge in $A$.) 
The solution of this local problem is obtained in 
three steps as in \S 3 of \cite{FCAP}.

\smallskip
\textit{Step 1:}
We denote by $f_0\colon X\to Y$ a map which is holomorphic 
on the set $A$ obtained in the earlier steps of the process. 
Write $F_0(x)=(x,f_0(x)) \in X\times Y$.
We find a small ball $0\in D\subset \C^p$ $(p=\dim Y)$
and a holomorphic map $f\colon A\times D\to Y$
such that $f(\cdotp,0)=f_0$, $f(x,t)=f_0(x)$ 
for all $x\in X_0$ and $t\in D$, and the partial 
differential $\di_t f(x,t) \colon T_t \C^r \to T_{f(x,t)} Y$ 
in the $t$ variable is surjective for $x$ in a 
neighborhood of $C=A\cap B$ and $t\in D$
(Lemma 3.2 in \cite{FCAP}). Such $f$ is found by 
choosing holomorphic vector fields $\xi_1,\ldots,\xi_p$
in a Stein neighborhood of $F_0(A)$ in $Z=X\times Y$
which are tanget the fibers $Y$, they vanish on the intersection 
of their domains with $X_0\times Y$ and span the 
vertical tangent bundle over a neighborhood 
of $F_0(C)$. The flow $\theta^j_t$ of $\xi_j$ exists 
for sufficiently small $t\in \C$, and the 
map 
$$
	F(x,t_1,\ldots,t_p) = (x,f(x,t))= 
	\theta^1_{t_1}\circ\cdots \circ 
	\theta^p_{t_p}\circ F_0(x)  \in X\times Y
$$
satisfies the stated requirements.

\textit{Step 2:} We approximate the holomorphic map 
$f$ from Step 1 uniformly on a neighborhood of 
the compact convex set $C\times D$ by a map 
$g$ which is holomorphic in an open neighborhood 
of $B\times D$. This is possible since $C$ is a compact 
convex set in $\C^n$ $(n=\dim X)$ with respect to some 
local holomorphic coordinates in a neighborhood of $B$ in $X$; 
hence $C\times D$ may be considered as a compact 
convex set in $\C^{n+p}$ and $g$ exists by the CAP property of $Y$. 
(In \cite{FCAP}, Definition 1.1, we used a more technical 
version of CAP which required approximability only on special 
compact convex sets with the purpose of making it easier to verify. 
However, by the main result of \cite{FCAP} this more restrictive 
definition implies the Oka property with approximation, hence in 
particular the approximability of holomorphic maps 
$K\to Y$ on any compact convex set $K\subset \C^m$ by entire 
maps $\C^m\to Y$.)

\textit{Step 3:} 
We `glue' $f$ and $g$ into a holomorphic map 
$f' \colon (A\cup B)\times D \to Y$ which approximates $f$ 
uniformly on $A\times D$ and agrees with it over $X_0$. 
(A small shrinking of the domain is necessary.) 
This is done by finding a fiberwise biholomorphic 
transition map $\gamma(x,t)$ which is close to $(x,t)\to t$
and satisfies $f(x,t)=g(x,\gamma(x,t))$ on $C\times D$, splitting
it in the form $\gamma_x=\gamma(x,\cdotp)=\beta_x\circ\alpha^{-1}_x$ 
$(x\in C)$ where $\alpha$ is holomorphic on $A\times D$ 
and $\beta$ is holomorphic on $B\times D$ (Lemma 2.1 in \cite{FCAP}),
and taking $f'(x,t)= f(x,\alpha(x,t))=g(x,\beta(x,t))$.
The map $f'_0=f'(\cdotp,0) \colon A\cup B \to Y$ 
is holomorphic on $A\cup B$, it approximates 
$f_0$ uniformly on $A$ and agrees with $f_0$ on~$X_0$
(the last observation follows from the construction of $f$). 
The construction also gives a homotopy with the required 
properties from $f_0$ to $f'_0$.
(We can also obtain $f'_0$ by applying Proposition 5.2 
in \cite{FP1} which is available in the parametric case; 
this can be used in the proof of the parametric analogue
of theorem \ref{t1}.)
\smallskip

This explains the {\em noncritical case} in which we are not
passing any nonsmooth boundary points of $bD_t$. The number
of steps needed to reach $D_{t_1}$ from $D_{t_0}$ 
does not depend on the partial solutions obtained in the 
intermediate steps. 

It remains to consider the critical values $t\in [0,1]$ for which 
$bD_t$ has a nonsmooth point in $\Omega$. 
The defining equation of $D_t\cap \Omega$ can be written as 
$$   
	D_t\cap \Omega = \{x\in \Omega \colon 
	h(x) \stackrel{def}{=} \frac{\tau(x)}{\tau(x)-\rho(x)} \le  t\}, 
$$
and the equation for critical points $dh=0$ is equivalent to 
$$
	(\tau -\rho)d\tau - \tau(d\tau -d\rho)= \tau d\rho-\rho d\tau =0.
$$
A generic choice of $\rho$ and $\tau$ insures that there are at most
finitely many solutions in $\Omega$ and no solution on $b\Omega$; 
furthermore, these solutions lie on distinct level sets of $h$. 
At each critical point of $h$ the complex Hessians satisfy  
$(\tau -\rho)^2 H_h = \tau H_\rho - \rho H_\tau$.
Since $\tau>0$ and $-\rho>0$ on $\Omega$, we conclude that 
$H_h>0$ at such points, i.e., $h$ is strongly plurisubharmonic 
near its critical point in $\Omega$. (See \cite{FActa} for more details.) 

A method for passing the isolated critical points of $h$ 
was explained in \S 6.2--\S 6.4 of \cite{FActa} for submersions $X\to\C^q$.
In \S  6 of \cite{FFourier} this method was adapted to holomorphic maps 
to any complex manifold $Y$; we include a brief description.  
Assume that $f\colon X\to Y$ is holomorphic on a 
certain sublevel set $A=D_{t_0}$ of $h$ just below a 
critical level $t_1$. We attach to $A$ a totally real handle 
$E \subset \Omega\bs X_0$ containing the critical 
(nonsmooth) point $p\in bD_{t_1}$ which describes the 
change of topology of $D_t$ as $t$ passes $t_1$. 
By theorem 3.2 in \cite{FFourier} the map $f|_{A\cup E}$ 
can be uniformly approximated by a map which is holomorphic 
in an open neighborhood of $A\cup E$ and agrees with $f$ on $X_0$. 
Finally, by the method explained in \S 6.4 of \cite{FActa} we 
can proceed by the noncritical case (applied with a different 
strongly plurisubharmonic function constructed especially for 
this purpose) to a sublevel set $D_{t_2}$ above 
the critical level ($t_2> t_1$). Now we can 
proceed to the next critical value of $h$.

This completes the proof of theorem \ref{t1} for maps $X\to Y$.
The same proof applies to sections of holomorphic fiber bundles
with fiber $Y$ over a Stein manifold $X$ because the use 
of CAP can be localized to small subsets in $X$ over which 
the bundle is trivial (compare with \cite{FCAP}).

%
%
%
%
\section{Proof of proposition \ref{p1}}
Let $Y$ be a complex manifold satisfying the conclusion
of theorem \ref{t1} for closed complex submanifolds in
Stein manifolds. Let $X$ be a Stein manifold, 
$K\subset X$ a compact $\cH(X)$-convex subset, 
$X_0$ a closed complex subvariety of $X$ 
and $f\colon X\to Y$ a continuous map which is holomorphic
in an open set $U_0\subset X$ containing $K\cup X_0$.

The set $K\cup X_0$ has a basis of open Stein neighborhoods 
in $X$ and hence we may assume (by shrinking $U_0$
if necessary) that $U_0$ is Stein. 
Choose a proper holomorphic embedding 
$\phi\colon U_0 \to\C^r$; its graph 
$\Sigma=\{(x,\phi(x))\colon x\in U_0\}$ is a 
closed complex submanifold of the Stein manifold $M=X\times \C^r$.
The map $g_0\colon M\to Y$, $g_0(x,z)=f_0(x)$ 
$(x\in X,\ z\in \C^r)$, is continuous on $M$ and holomorphic 
on the open set $U_0 \times \C^r\subset M$ containing $\Sigma$. 

By the assumption on $Y$ there is a homotopy 
$g_t\colon M\to Y$ $(t\in [0,1])$ which is fixed 
on the submanifold $\Sigma$ and such that $g_1$ is holomorphic 
on $M$. For later purposes we also need that $g_t$ be 
holomorphic in an open neighborhood $V \subset X\times Y$
(independent of $t\in [0,1]$) of $\Sigma$. 
Since both $g_0$ and $g_1$ are holomorphic near $\Sigma$, 
this can be accomplished (without changing $g_0$ and $g_1$)
by a simple modification of the homotopy $g_t$, using a 
strong holomorphic deformation retraction of a neighborhood 
of $\Sigma$ in $M$ onto $\Sigma$ (such exists by combining 
Siu's theorem \cite{Si}, \cite{De} with the Docquier-Grauert 
theorem \cite{DG}; see Corollary 1 in \cite{Si}).

To complete the proof we shall need the following 
(known) result which follows from Cartan's theorems A and B
(see lemma 8.1 in \cite{FP2} for a reduction 
to the Oka-Weil approximation theorem).

\begin{lemma} 
\label{l1}
Let $X$ be a Stein manifold, $L\subset X$ 
a compact $\cH(X)$-convex subset and $X_0\subset X$ a 
closed complex subvariety. Let $\phi$ be a holomorphic function
in an open set containing $L\cup X_0$. For every $s\in \N$ 
and $\eta >0$ there exists $\varphi\in \cH(X)$ such that $\varphi -\phi$ 
vanished to order $s$ on $X_0$ and 
$\sup_{x\in L} |\varphi(x)-\phi(x)| <\eta$.
\end{lemma}

We apply this lemma to the embedding $\phi\colon U\to \C^r$
and a compact set $L\subset U_0$ with $K\subset {\rm Int}\,L$.
Denote the resulting map by $\varphi\colon X\to\C^r$. 
Consider the homotopy $f_t\colon X\to Y$ 
defined by $f_t(x)= g_t(x,\varphi(x))$ for $x\in X$ and $t\in [0,1]$.
At $t=0$ we get the initial map $f_0(x)=g_0(x,\varphi(x))$, and 
at $t=1$ we get a map $f_1(x)=g_1(x,\varphi(x))$ which is holomorphic
on $X$. There is a small open neighborhood $U\subset U_0$
of $K\cup X_0$ such that $(x,\varphi(x))\in V$ for $x\in U$;
since $g_t$ is holomorphic on $V$, we see that $f_t$ is holomorphic 
on $U$. It is easily verified that for every $t\in [0,1]$,
$f_t$ agrees with $f_0$ to order $s$ along $X_0$, and it 
approximates $f_0$ uniformly on $K$ as well as desired provided 
that the number $\eta$ in lemma\ref{l1} is chosen sufficiently small.
Thus $f_t$ satisfies the conclusion of proposition \ref{p1}.

All steps go through in the parametric case and give
the corresponding parametric analogue of proposition \ref{p1}.

\begin{remark} 
\label{p1sections}
The analogoue of proposition \ref{p1} holds for sections 
of fiber bundles with fiber $Y$ over Stein manifolds:
interpolation implies jet interpolation and approximation. 
To see this, assume that $\pi\colon E\to X$ is such a bundle, 
$X_0$ is a closed complex subvariety of $X$ 
and $K$ is a compact $\cH(X)$-convex subset of $X$.
Let $f_0\colon X\to E$ be a continuous section
which is holomorphic in an open Stein neighborhood
$U\subset X$ of $K\cup X_0$. As before we embed $U$ to $\C^r$ 
and denote the graph of the embedding by 
$\Sigma\subset M=X\times \C^r$.  Let $p^*E\to M$ 
be the pull-back of $E$ by the projection 
$p\colon M\to X$; this bundle has the same 
fiber $Y$ and there is a natural map $\theta\colon p^* E\to E$ 
covering $p$ which is isomorphism on the fibers. 
The section $f_0$ pulls back to a section 
$h_0\colon M\to p^* E$ satisfying $\theta(h_0(x,y))=f_0(x)$.
Note that $h_0$ is holomorphic on $U'=p^{-1}(U)=U\times\C^r$. 
By the hypothesis there is a homotopy $h_t\colon M\to p^*E$ 
which is fixed on $\Sigma$ such that $h_1$ is holomorphic.
Choosing $\varphi \colon X\to\C^r$ as before we obtain 
a homotopy $f_t(x)= \theta (h_t(x,\varphi(x))) \in E_x$ 
$(x\in X)$  satisfying the desired conclusion. 
\end{remark}

\begin{remark}  L\'arusson proved that for a Stein manifold,  
the basic Oka property with interpolation implies ellipticity, 
i.e., existence of a dominating spray (Theorem 2 in \cite{L}). 
This also follows from proposition \ref{p1} and 
the  known result that ellipticity of a Stein manifold 
is implied by second order jet interpolation 
(\cite{GOka}; proposition 1.2 in \cite{FP3}). 
\end{remark}

\section{Equivalences between the parametric Oka properties}
In this section we discuss the parametric 
analogue of corollary \ref{c1}.

Let $P$ be a compact Hausdorff space 
(the parameter space) and $P_0\subset P$ a closed subset of $P$ 
which is a strong deformation retract of some open 
neighborhood of $P_0$ in $P$. In most applications $P$ is a 
polyhedron and $P_0$ a subpolyhedron; an important
special case is $P_0=\{0,1\} \subset [0,1]=P$. 

Given a Stein manifold $X$ and a compact $\cH(X)$-convex 
subset $K$ in $X$ we consider continuous maps $f\colon X\times P\to Y$ 
such that for every $p\in P$ the map $f^p=f(\cdotp,p)\colon X\to Y$ 
is holomorphic in an open neighborhood of $K$ in $X$ 
(independent of $p\in P$), and for every $p\in P_0$ 
the map $f^p$ is holomorphic on $X$. We say that $Y$ satisfies the 
{\em parametric Oka property with approximation} 
if for every such data $(X,K,P,P_0,f)$ there is a homotopy 
$f_t\colon X\times P\to Y$ $(t\in [0,1])$ consisting
of maps satisfying the same properties as $f_0=f$ such that

\begin{itemize}
\item[(i)]  
the homotopy is fixed on $P_0$: $f_t^p=f_0^p$ when $p\in P_0$ 
and $t\in [0,1]$, 
\item[(ii)] 
$f_t$ approximates $f_0$ uniformly on $K\times P$ for all $t\in [0,1]$, and
\item[(iii)] 
$f_1^p \colon X\to Y$ is holomorphic for every $p\in P$.
\end{itemize}

We say that $Y$ satisfies the 
{\em parametric convex approximation property} (PCAP) 
if the above holds for any compact convex set 
$K\subset X=\C^n$, $n\in \N$. 

Similarly, $Y$ satisfies the 
{\em parametric Oka property with interpolation}
if for any closed complex subvariety $X_0$ in a 
Stein manifold $X$ and continuous map 
$f\colon X\times P\to Y$ such that $f^p|_{X_0}$ is 
holomorphic for all $p\in P$ and $f^p$ is holomorphic
on $X$ for $p\in P_0$ there is a homotopy $f_t\colon X\times P\to Y$ 
$(t\in [0,1])$, with $f_0=f$, which satisfies 
(i) and (iii) above together with the interpolation condition

\begin{itemize}
\item[(ii')]  $f_t^p(x)=f_0^p(x)$ for $x\in X_0$, $p\in P$ and $t\in [0,1]$.
\end{itemize}

Similarly one introduces the {\em parametric Oka property 
with jet interpolation} on subvarieties $X_0\subset X$.
Combining it with approximation on $\cH(X)$-convex sets one obtains 
{\em the parametric Oka property with jet interpolation and approximation}
which coincides with Gromov's ${\rm Ell}_\infty$ property (\cite{GOka}, \S 3.1). 
By Theorem 1.5 in \cite{FP2} and Theorem 1.4 in \cite{FP3} all these properties 
are implied by ellipticity of $Y$ (the existence of a dominating spray), 
and also by subellipticity of $Y$ \cite{Fsubell}.
In analogy to the non-parametric case (corollary \ref{c1})
these ostensibly different properties are equivalent.

\begin{theorem}
\label{parametric}
The following properties of a complex manifold are equivalent:
\begin{enumerate}
\item[(a)] the parametric convex approximation property {\rm (PCAP)};
\item[(b)] the parametric Oka property with approximation;
\item[(c)] the parametric Oka property with interpolation;
\item[(d)] the parametric Oka property with jet interpolation and 
approximation, i.e., Gromov's ${\rm Ell}_\infty$ property 
{\rm (\cite{GOka},  \S 3.1)}.
\end{enumerate}
\end{theorem}

\begin{proof}
The logic of proof is the same as in corollary \ref{c1}:
(a)$\Leftrightarrow$(b) is Theorem 5.1 in \cite{FCAP},
(a)$\Rightarrow$(c) is the parametric version of theorem 
\ref{t1} in this paper, and (c)$\Rightarrow$(d) is the 
parametric version of proposition \ref{p1}.
The remaining implications are consequences of definitions
and of Proposition 1.2 in \cite{FP3} concerning the holomorphic 
extendability of holomorphic maps from a subvariety 
$X_0\subset X$ to an open neighborhood of $X_0$ in $X$. 
The implication (c)$\Rightarrow$(b) 
was first proved by L\'arusson \cite{L}.

To obtain these implications in the parametric case it suffices 
to observe that theorem \ref{t1} and proposition \ref{p1} 
remain valid in this more general situation as we have already 
indicated at the appropriate places in their proofs.
The reader who may be interested in further details 
should consult the proof of the parametric Oka property 
with approximation for elliptic manifolds in \cite{FP1} 
(Theorems 1.4 and 5.5).

The following observation might be helpful
concerning the use of theorem \ref{Steinnbds} 
(the existence of Stein neighborhoods)
in the proof of  (a)$\Rightarrow$(c) in theorem \ref{parametric}. 
In the proof of proposition \ref{pmain} (\S 4) we needed 
a Stein neighborhood $U\subset X\times Y$ of the graph of 
a holomorphic map $K\cup X_0\to Y$.
It may appear that one  needs a stronger theorem in the 
parametric case (a family of Stein neighborhoods 
depending continuously on the parameter $p\in P$). 
This is unnecessary because the graphs of maps $f^p$ 
over the set $K'= K\cup (L\cap X_0)$ in the proof of 
proposition \ref{pmain} belong to the same Stein set 
$U$ for $p\in P$ close to some initial point $p_0$ for 
which $U$ was chosen. Finitely many such Stein open sets 
in $X\times Y$ cover the graph of $f\colon X\times P\to Y$ over the
compact set $K'\times P$; the solutions  in each of these 
Stein neighborhoods (for an open set of parameters in $P$) 
are patched together by a continuous partition 
of unity in the parameter space (see the proof of Theorem 5.5 
in \cite{FP1}). 
\end{proof}

A continuous map satisfying the homotopy lifting property
is called a {\em Serre fibration} \cite{W}, p.\ 8.   
A holomorphic submersion $\pi\colon Y\to Y_0$ is said to
be {\em subelliptic} if every point $y_0\in Y_0$
has an open neighborhood $U\subset Y_0$ such that
$\pi\colon \pi^{-1}(U)\to U$ admits a finite 
fiber-dominating family of sprays (\cite{Fsubell}, p.\ 529); 
for a holomorphic fiber bundle this is equivalent to 
subellipticity of the fiber. Theorem \ref{parametric}
together with the parametric version of Theorem 1.3 
in \cite{FCAP} (see the discussion following 
Theorem 5.1 in \cite{FCAP}) we obtain the following.

\begin{corollary} 
\label{c2}
Suppose that $\pi \colon Y\to Y_0$ is a subelliptic submersion 
which is also a Serre fibration. If $Y_0$ satisfies any 
(and hence all) of the properties in theorem \ref{parametric} 
then so does $Y$. Conversely, if $\pi$ is surjective 
and $Y$ satisfies any of these properties with a contractible 
parameter space $P$ then so does $Y_0$. 
\end{corollary}

\smallskip
\noindent \textbf{A concluding remark.}
Our results suggest that CAP (and its parametric analogue) 
is the most natural Oka-type property to be studied further since
it is the simplest to verify, yet equivalent to all other 
Oka properties. Indeed CAP is just the localization of 
the {\it Oka property with approximation} to maps 
from compact convex sets in Eulidean spaces 
(see \cite{Fflex} for this point of view).
CAP easily follows from ellipticity or subellipticity,
and is a natural opposite property to Kobayashi-Eisenman-Brody hyperbolicity 
\cite{Br}, \cite{E}, \cite{Kb}. It is also related to dominability 
by Euclidean spaces, a property that has been extensively studied 
\cite{BL}, \cite{CG},  \cite{KO}, \cite{Kd}.
It would be interesting to know how big (if any)
is the gap between CAP and (sub)ellipticity.
Several related questions are mentioned in 
\cite{Fsurvey} and  \cite{Fflex}.

\medskip \textit{Acknowledgement.} I thank
F.\ L\'arusson for useful discussions; in particular, 
the proof of proposition \ref{p1} is motivated 
by Theorem 1 in his recent preprint \cite{L}.
I also thank B.\ Drinovec-Drnov\v sek from whom I have 
learned the geometric scheme (pushing out by a convex
combination of two plurisubharmonic functions) 
used in the proof of proposition \ref{pmain} 
and in \cite{FActa}, p.\ 179, and J.\ Prezelj-Perman
for pointing out an inaccuracy in an earlier version
of the paper. Finally, I thank M.\ Col\c toiu for pointing 
out the references \cite{Co}, \cite{CM} and \cite{Pe} in 
connection with theorem \ref{Steinnbds}, and 
C.\ Laurent Thi\'ebaut for the translation 
of the abstract in French.

\bibliographystyle{amsplain}

\end{document}